\documentclass{amsart}

\usepackage{amsmath} 
\usepackage{amssymb}

\newtheorem{theorem}{Theorem}[section] 
\newtheorem{claim}{Claim}[theorem]
\newtheorem{lemma}[theorem]{Lemma} 
\newtheorem{proposition}[theorem]{Proposition} 
\newtheorem{corollary}[theorem]{Corollary} 

\theoremstyle{definition}
\newtheorem{definition}[theorem]{Definition}

\newtheorem{problem}[theorem]{Problem}

\theoremstyle{remark}
\newtheorem{remark}[theorem]{Remark}

\numberwithin{equation}{section}
\newcommand{\forces}{\Vdash}

\newcommand{\lbv}{[\![} 
\newcommand{\rbv}{]\!]}
\newcommand{\bV}{{\bf V}} 
\newcommand{\lesdot}{\mathrel{\mathord{<}\!\!\raise 
0.8 pt\hbox{$\scriptstyle\circ$}}} 
\newcommand{\bO}{\emptyset}

\newcommand{\gb}{{\mathfrak b}}

\newcommand{\can}{2^{\textstyle \omega}} 
\newcommand{\fs}{2^{\textstyle <\!\omega}} 
\newcommand{\baire}{\omega^{\textstyle \omega}} 
\newcommand{\iso}{[\omega]^{\textstyle \omega}} 
 
\newcommand{\fseo}{\omega^{\textstyle <\!\omega}} 

\newcommand{\conc}{{}^\frown\!}
\newcommand{\lh}{{\rm lh}\/}
\newcommand{\rest}{{\restriction}}
\newcommand{\mrot}{{\rm root}\/} 
\newcommand{\suc}{{\rm succ}} 
\newcommand{\spliting}{{\rm split}}

\newcommand{\cls}{{\rm cl}^{\rm sw}}
\newcommand{\mul}{{\mu^{\rm Leb}}}
 
\newcommand{\dn}{{\rm dn}}
\newcommand{\up}{{\rm up}}

\newcommand{\borel}{{\bf Borel}}
\newcommand{\lev}{{\rm lev}}
\newcommand{\UM}{{\mathbb U}{\mathbb M}}

\newcommand{\bqsc}{{\bbQ^{\cT^{\rm sc}}}}
\newcommand{\rsc}{{\rm sc}}
\newcommand{\bqT}{{\bbQ^{{\mathcal T}}}}
\newcommand{\rk}{{\rm rk}}

\newcommand{\cA}{{\mathcal A}}
\newcommand{\bbA}{{\mathbb A}}

\newcommand{\cB}{{\mathcal B}}
\newcommand{\BB}{{\mathbb B}}

\newcommand{\bbC}{{\mathbb C}}

\newcommand{\bbD}{{\mathbb D}}

\newcommand{\bH}{{\bf H}}

\newcommand{\cF}{{\mathcal F}}
\newcommand{\cG}{{\mathcal G}}
\newcommand{\cI}{{\mathcal I}}

\newcommand{\bbL}{{\mathbb L}}

\newcommand{\bbP}{{\mathbb P}}

\newcommand{\gp}{{\mathfrak p}}

\newcommand{\bbQ}{{\mathbb Q}}

\newcommand{\mbR}{{\mathbb R}}

\newcommand{\cT}{{\mathcal T}}
\newcommand{\cU}{{\mathcal U}}

\newcount\skewfactor
\def\mathunderaccent#1#2 {\let\theaccent#1\skewfactor#2
\mathpalette\putaccentunder}
\def\putaccentunder#1#2{\oalign{$#1#2$\crcr\hidewidth
\vbox to.2ex{\hbox{$#1\skew\skewfactor\theaccent{}$}\vss}\hidewidth}}
\def\name{\mathunderaccent\tilde-3 }

\begin{document}

\title{How much sweetness is there in the universe?}

\author{Andrzej Ros{\l}anowski}
\address{Department of Mathematics\\
 University of Nebraska at Omaha\\
 Omaha, NE 68182-0243, USA}
\email{roslanow@member.ams.org}
\urladdr{http://www.unomaha.edu/$\sim$aroslano}

\author{Saharon Shelah}
\address{Einstein Institute of Mathematics\\
Edmond J. Safra Campus, Givat Ram\\
The Hebrew University of Jerusalem\\
Jerusalem, 91904, Israel\\
 and  Department of Mathematics\\
 Rutgers University\\
 New Brunswick, NJ 08854, USA}
\email{shelah@math.huji.ac.il}
\urladdr{http://www.math.rutgers.edu/$\sim$shelah}
\thanks{Both authors acknowledge support from the United States-Israel
Binational Science Foundation (Grant no. 2002323). This is publication 856
of the second author.}

\date{September 2004}

\begin{abstract}
We continue investigations of forcing notions with strong ccc properties
introducing new methods of building sweet forcing notions. We also show that
quotients of topologically sweet forcing notions over Cohen reals are
topologically sweet while the quotients over random reals do not have to be
such.
\end{abstract}

\maketitle

\section{Introduction}
One of the main ingredients of the construction of the model for {\em all
projective sets of reals have the Baire property\/} presented in Shelah
\cite[\S 7]{Sh:176} was a strong ccc property of forcing notions called {\em
sweetness}. This property is preserved in amalgamations and also in
compositions with the Hechler forcing notion $\bbD$ and the Universal Meager
forcing $\UM$ (see \cite[\S 7]{Sh:176}; a full explanation of how this is
applied can be found in \cite{JuRo}). Stern \cite{St85} considered a
slightly weaker property, {\em topological sweetness}, which is also
preserved in amalgamations and compositions with $\bbD$ and $\UM$. We
further investigated the sweet properties of forcing notions in \cite[\S
4]{RoSh:672}, where we introduced a new property called {\em iterable
sweetness\/} and we showed how one can build sweet forcing notions. New
examples of iterably sweet forcing notions can be used in constructions like 
\cite[\S 7]{Sh:176}, \cite{Sh:F380}, but it could be that there is no need
for this --- the old forcing notions could be adding generic objects for all
of them. In \cite{RoSh:845} we proved that this is exactly what happens with
the natural examples of sweet forcing notions determined by the universality
parameters as in \cite[\S 2.3]{RoSh:672}: a sequence  Cohen real ---
dominating real --- Cohen real produces generic filters for many of them. 

In the present paper we show that sweetness is not so rare after all and we
give more constructions of sweet forcing notions. In the first section we
present a new method of building sweet forcing notions and we give our first
example:  a forcing notion $\bqsc$ associated with scattered subtrees of
$\fs$. We do not know if the iterations of ``old'' forcing notions add
generic objects for $\bqsc$, but in Proposition \ref{nogen} we present an
indication that this does not happen. In the second section we use our
method to introduce two large families of sweet forcing notions, in some
sense generalizing the known examples from \cite{RoSh:672}. This time we
manage to show that {\em some\/} of our forcing notions are really new by
showing that we have too many different examples (in Theorems
\ref{secondnew}, \ref{newthm}). 

In the last section of the paper we investigate the preservation of
topological sweetness under some operations. We note that a complete
subforcing of a topologically sweet separable partial order is equivalent to
a topologically sweet forcing (in Proposition \ref{tosub}). We also show that
the quotient of a topologically sweet forcing notion by a Cohen subforcing
is topologically sweet (Theorem \ref{overc}), but quotients by random real
do not have to be topologically sweet (Corollary \ref{overran}).  

\subsection{Notation}
Our notation is rather standard and compatible with that of classical
textbooks (like Jech \cite{J} or Bartoszy\'nski and Judah \cite{BaJu95}). In
forcing we keep the older convention that {\em a stronger condition is the
larger one}. Our main conventions are listed below.  

\begin{enumerate}
\item For a forcing notion $\bbP$, $\Gamma_\bbP$ stands for the canonical
$\bbP$--name for the generic filter in $\bbP$. With this one exception, all
$\bbP$--names for objects in the extension via $\bbP$ will be denoted with
a tilde below (e.g., $\name{\tau}$, $\name{X}$). 

\noindent The weakest element of $\bbP$ will be denoted by $\bO_\bbP$ (and
we will always assume that there is one and that there is no other condition 
equivalent to it).  

\item The complete Boolean algebra determined by a forcing notion $\bbP$ is 
denoted by ${\bf BA}(\bbP)$. For a complete Boolean algebra $\BB$, $\BB^+$
is $\BB\setminus\{{\bf 0}_{\BB}\}$ treated as a forcing notion (so the
order is the reverse Boolean order). Also, for a formula $\varphi$, the
Boolean value (with respect to $\BB$) of $\varphi$ will be denoted by
$\lbv\varphi\rbv_\BB$.  

\item Ordinal numbers will be denoted be the lower case initial letters of
the Greek alphabet ($\alpha,\beta,\gamma,\delta\ldots$) and also by $i,j$
(with possible sub- and superscripts).\\ 
Cardinal numbers will be called $\kappa,\lambda,\mu$. 

\item For two sequences $\eta,\nu$ we write $\nu\vartriangleleft\eta$
whenever $\nu$ is a proper initial segment of $\eta$, and $\nu
\trianglelefteq\eta$ when either $\nu\vartriangleleft\eta$ or $\nu=\eta$. 
The length of a sequence $\eta$ is denoted by $\lh(\eta)$.

\item The quantifier $(\exists^\infty n)$ is an abbreviation for $(\forall
  m\in\omega) (\exists n>m)$. 

\item The Cantor space $\can$ and the Baire space $\baire$ are the spaces of
all functions from $\omega$ to $2$, $\omega$, respectively, equipped with
the natural (Polish) topology. 
\end{enumerate}

\subsection{Background on sweetness}
Let us recall basic definitions related to sweet forcing notions.

\begin{definition}
[Shelah {\cite[Def.~7.2]{Sh:176}}] 
\label{sweet}
A pair $(\bbP,\bar{E})$ is {\em model of sweetness\/} whenever: 
\begin{enumerate}
\item[(i)]    $\bbP$ is a forcing notion, 
\item[(ii)]   $\bar{E}=\langle E_n:n<\omega\rangle$, each $E_n$ is an
equivalence relation on $\bbP$ such that $\bbP/E_n$ is countable,
\item[(iii)]  equivalence classes of each $E_n$ are $\leq_{\bbP}$--directed, 
$E_{n+1}\subseteq E_n$,
\item[(iv)]   if $\{p_i:i\leq\omega\}\subseteq\bbP$, $p_i\; E_i\; p_\omega$
(for $i\in\omega$), then
\[(\forall n\in\omega)(\exists q\geq p_\omega)(q\; E_n\; p_\omega\ \&\
(\forall i\geq n)(p_i\leq q)),\]
\item[(v)]    if $p,q\in\bbP$, $p\leq q$ and $n\in\omega$, then there is
$k\in\omega$ such that
\[(\forall p'\in [p]_{E_k})(\exists q'\in [q]_{E_n})(p'\leq q').\]
\end{enumerate}
If there is a model of sweetness based on $\bbP$, then we say that $\bbP$ is
{\em sweet}. 
\end{definition}

\begin{definition}
[Stern {\cite[Def.~1.2]{St85}}]
\label{topsweet}
{\em A model of topological sweetness\/} is a pair $M=(\bbP,\cB)$ such
that $\bbP=(\bbP,\leq)$ is a forcing notion, $\cB$ is a countable basis
of a topology $\tau$ on $\bbP$ and  
\begin{enumerate}
\item[(i)]   $\bO_\bbP$ is an isolated point in $\tau$,
\item[(ii)]  if a sequence $\langle p_n:n<\omega\rangle\subseteq\bbP$ is
$\tau$--converging to $p\in\bbP$, $q\geq p$ and $W$ is a
$\tau$--neighbourhood of $q$, then there is a condition $r\in\bbP$ such that  
\begin{enumerate}
\item[(a)] $r\in W$, $r\geq q$,
\item[(b)] the set $\{n\in\omega: p_n\leq r\}$ is infinite.
\end{enumerate}
\end{enumerate}
If there is a model of topological sweetness $(\bbP,\cB)$, then the forcing
notion $\bbP$ is {\em topologically sweet}.
\end{definition}

\begin{lemma}
[See {\cite[Lemma 4.2.3]{RoSh:672}}]
\label{bastoplem}
Assume that $(\bbP,\cB)$ is a model of topological sweetness.
\begin{enumerate}
\item If $p,q\in\bbP$, $p\leq q$ and $q\in U\in\cB$, then there is an open
neighbourhood $V$ of $p$ such that
\[(\forall r\in V)(\exists r'\in U)(r\leq r').\]
\item If $m\in\omega$, $p\in U\in \cB$, then there is an open neighbourhood 
$V$ of $p$ such that any $p_0,\ldots,p_m\in V$ have a common upper bound in
$U$.
\end{enumerate}
\end{lemma}

\begin{definition}
[See {\cite[Def. 4.2.1]{RoSh:672}}]
\label{itersweet}
Let $\cB$ be a countable basis of a topology on a forcing notion $\bbQ$. We
say that $(\bbQ,\cB)$ is {\em a model of iterable sweetness\/} if 
\begin{enumerate}
\item[(i)] $\cB$ is closed under finite intersections,
\item[(ii)] each $U\in \cB$ is directed and $p\leq q\in U\ \Rightarrow\
p\in U$,  
\item[(iii)] if $\langle p_n:n\leq \omega\rangle\subseteq U$ and the
sequence $\langle p_n:n<\omega\rangle$ converges to $p_\omega$ (in the
topology generated by $\cB$), then there is a condition $p\in U$ such that
$(\forall n\leq\omega)(p_n\leq p)$. 
\end{enumerate}
\end{definition}

\begin{proposition}
[See {\cite[Proposition 4.2.2]{RoSh:672}}]
If $\bbP$ is a sweet forcing notion in which any two compatible conditions
have a least upper bound, then $\bbP$ is iterably sweet.
\end{proposition}

\section{sw--closed families and scattered trees}
In this section we present a new method of building sweet forcing
notions. This method is, essentially, a generalization of that determined by
the universality parameters of \cite[\S 2.3]{RoSh:672}.

\begin{definition}
\label{tree}
\begin{enumerate}
\item A {\em tree} is a family $T$ of finite sequences such that for some
$\mrot(T)\in T$ we have
\[(\forall\nu\in T)(\mrot(T)\trianglelefteq \nu)\quad\mbox{ and }\quad
\mrot(T)\trianglelefteq\nu\trianglelefteq\eta\in T\ \Rightarrow\ \nu\in T.\]
\item If $\eta$ is a node in the tree $T$ then 
\[\begin{array}{lcl}
\suc_T(\eta)&=&\{\nu\in T: \eta\vartriangleleft\nu\ \&\ \lh(\nu)=\lh(\eta)+1
\}\ \mbox{ and}\\
T^{[\eta]}&=&\{\nu\in T:\eta\trianglelefteq\nu\}.
  \end{array}\]
\item For a tree $T$, the family of all $\omega$--branches through $T$ is
denoted by $[T]$, and we let
\[\max(T)\stackrel{\rm def}{=}\{\nu\in T:\mbox{ there is no }\rho\in
T\mbox{ such that }\nu\vartriangleleft\rho\}\]   
and 
\[\spliting(T)\stackrel{\rm def}{=}\{\nu\in T:|\suc_T(\nu)|\geq 2\}.\]   
\item A tree $T$ is {\em normal\/} if $\max(T)=\emptyset$ and $\mrot(T)=\langle
\rangle$. 
\end{enumerate}
\end{definition}

\begin{definition}
\label{closefam}
Suppose that $\cT$ is a family of normal subtrees of $\fseo$. We say that
$\cT$ is {\em sw--closed\/} whenever  
\begin{enumerate}
\item if $T_1\in\cT$, $T_2\subseteq T_1$ and $T_2$ is a normal tree, then
  $T_2\in\cT$, 
\item if $T_1,T_2\in \cT$, then $T_1\cup T_2\in\cT$, and 
\item if $\langle T_n:n\leq\omega\rangle\subseteq\cT$ is such that
$\big(\forall n<\omega\big)\big(T_\omega\cap \omega^{\textstyle {\leq}n}=
T_n\cap\omega^{\textstyle {\leq} n}\big)$,\\
then $\bigcup\limits_{n\leq\omega} T_n\in \cT$.    
\end{enumerate}
\end{definition}

\begin{definition}
\label{genscat}
For a family $\cT$ of normal subtrees of $\fseo$ we define a forcing notion
$\bqT$ as follows.\\ 
{\bf A condition in $\bqT$} is a pair $p=(N^p,T^p)$ such that $N^p<\omega$
and $T^p\in\cT$.\\
{\bf The order $\leq_{\bqT}$ of $\bqT$} is given by 

$p\leq_{\bqT} q$\qquad if and only if 

$N^p\leq N^q$, $T^p\subseteq T^q$ and $T^q\cap \omega^{\textstyle N^p}=
T^p\cap \omega^{\textstyle N^p}$.
\end{definition}

The relation between the forcing $\bqT$ and the family $\cT$ is similar to
that in the case of the Universal Meager forcing notion $\UM$ and nowhere
dense subtrees of $\fs$. Note that $\bqT$ does not have to be ccc in
general, however in many natural cases it is.

\begin{proposition}
\label{swimsw}
Assume that $\cT$ is an sw--closed family of normal subtrees of $\fseo$ such 
that every $T\in\cT$ is finitely branching. Then $\bqT$ is a sweet forcing
notion in which any two compatible conditions have a least upper bound (and 
consequently $\bqT$ is iterably sweet).
\end{proposition}

\begin{proof}
One easily verifies that $\bqT$ is indeed a forcing notion and that any two 
compatible conditions in $\bqT$ have a least upper bound. 

\noindent For an integer $n<\omega$ let $E_n$ be a binary relation on
$\bqT$ defined by

$q\; E_n\; p$\qquad if and only if 

$N^q=N^p$ and $T^q\cap \omega^{\textstyle\leq N^q+n}=T^p\cap
  \omega^{\textstyle\leq N^q+n}$,

\noindent and let $\bar{E}=\langle E_n:n<\omega\rangle$. We claim that
$(\bqT,\bar{E})$ is a model of sweetness. Conditions \ref{sweet}(i--iii)
should be clear. To verify  \ref{sweet}(iv) suppose that $p_i\in\bqT$ for
$n\leq i\leq\omega$ are such that $p_i\; E_i\; p_\omega$ (for
$i<\omega$). Thus, for $n\leq i<\omega$, $N^{p_i}=N^{p_\omega}$ and
\[T^{p_i}\cap \omega^{\textstyle\leq N^{p_i}+i}=T^{p_\omega}\cap
\omega^{\textstyle\leq N^{p_\omega}+i}.\]
Put $N=N^{p_\omega}$ and $T=\bigcup\{T^{p_i}:n\leq i\leq\omega\}$. It
follows from \ref{closefam}(3) that $T\in\cT$, and plainly
$q=(N,T)\in\bqT$, $q\; E_n\; p_\omega$ and $(\forall i\geq n)(p_i\leq q)$,
finishing justification of \ref{sweet}(iv).  

Finally, to check \ref{sweet}(v) suppose that $p,q\in\bqT$, $p\leq q$ and
$n<\omega$. Let $k=N^q+n$. It should be clear that $(\forall p'\in
[p]_{E_k})(\exists q'\in [q]_{E_n})(p'\leq q')$.
\end{proof}

Now we are going to present our first example of an sw--closed family:
the family of scattered subtrees of $\fs$. 

\begin{definition}
\begin{enumerate}
\item For a closed set $A\subseteq\can$, let $\rk(A)$ be the
  Cantor--Bendixson rank of $A$, that is 
\[\rk(A)=\min\{\alpha<\omega_1:A^\alpha=A^{\alpha+1}\},\]
where $A^\alpha$ denotes the $\alpha^{\rm th}$ Cantor--Bendixson derivative 
of $A$. 
\item We say that a tree $T\subseteq\fs$ is {\em scattered} if it is normal
and $[T]$ is countable. The family of all scattered subtrees of $\fs$ will
be denoted by $\cT^\rsc$.
\item For a scattered tree $T\subseteq\fs$, let $g^T:[T]\longrightarrow
\rk([T])$ and $h^T:[T]\longrightarrow\omega$ be such that for each $\eta\in
[T]$ we have 
\[g^T(\eta)=\min\{\alpha<\rk(T):\eta\notin [T]^{\alpha+1}\}\]
and
\[h^T(\eta)=\min\{m<\omega:\big(\forall\nu\in [T]\big)\big(\nu\rest
m=\eta\rest m\ \Rightarrow\ (\eta=\nu\ \vee\ g^T(\nu)<g^T(\eta)\big)\}.\]
\end{enumerate}
\end{definition}

\begin{proposition}
Let $T\subseteq\fs$ be a normal tree. Then $T$ is scattered if and only if
there is a mapping $\varphi:T\longrightarrow\omega_1$ such that 
\begin{enumerate}
\item[$(\circledast)^0_{\varphi,T}$] 
\quad $\big(\forall \eta,\nu\in T\big)\big(\nu\vartriangleleft \eta\
\Rightarrow\ \varphi(\nu)\geq \varphi(\eta)\big)$, and 
\item[$(\circledast)^1_{\varphi,T}$] 
\quad $\big(\forall \eta\in \spliting(T)\big)\big(\varphi(\eta\conc\langle 
0\rangle)<\varphi(\eta)\ \vee\ \varphi(\eta\conc\langle 1\rangle)<\varphi
(\eta)\big)$.
\end{enumerate}
\end{proposition}

\begin{proof}
It should be clear that if there is a function $\varphi:T\longrightarrow
\omega_1$ such that $(\circledast)^0_{\varphi,T}+(\circledast)^1_{\varphi,
T}$ holds true, then the tree $T$ contains no perfect subtree and hence $T$
is scattered.  

We will show the converse implication by induction on $\rk(T)$. 

Suppose that $T$ is a scattered tree. Choose $\{\eta_\ell:\ell<n\}\subseteq
[T]$, $n<\omega$, such that $F\stackrel{\rm def}{=}\{\eta_\ell\rest
h^T(\eta_\ell):\ell<n\}$ is a front of $T$ and let 
\[A\stackrel{\rm def}{=}\big\{\rho\in T:\big(\exists\ell<n\big)\big(h^T(
\eta_\ell)<\lh(\rho)\ \&\ \rho\rest(\lh(\rho)-1)\vartriangleleft\eta_\ell\
\&\ \rho\ntriangleleft \eta_\ell\big)\big\}.\] 
Note that if $\ell<n$, $\nu\in [T]\setminus\{\eta_\ell\}$ and $\nu\rest h^T(
\eta_\ell)=\eta_\ell\rest h^T(\eta_\ell)$, then $g^T(\nu)<g^T(\eta_\ell)$. 
Hence $\big(\forall\rho\in A\big)\big(\rk(T^{[\rho]})<\rk(T)\big)$, 
so by the inductive hypothesis for each $\nu\in A$ we may choose
$\varphi_\nu:T^{[\nu]}\longrightarrow\omega_1$ such that $(\circledast)^0_{
\varphi_\nu,T^{[\nu]}}+(\circledast)^1_{\varphi_\nu,T^{[\nu]}}$ holds
true. Put $\alpha^*=\sup\{\varphi_\nu(\nu):\nu\in A\}<\omega_1$, $k^*=\max\{ 
h^T(\eta_\ell):\ell<n\}+1$  and let $\varphi:T\longrightarrow\omega_1$ be
defined by  
\[\varphi(\eta)=\left\{\begin{array}{ll}
\alpha^*+k^*-\lh(\eta)&\mbox{ if no initial segment of $\eta$ belongs to
  $F$, and}\\   
\alpha^*+1&\mbox{ if an initial segment of $\eta$ belongs to $F$}\\
          &\mbox{ but no initial segment of $\eta$ belongs to $A$, and}\\   
\varphi_\nu(\eta)&\mbox { if }\nu\in A\mbox{ and }\nu\trianglelefteq\eta.\\
		       \end{array}\right.\]
One easily verifies that the function $\varphi$ (is well defined and)
satisfies $(\circledast)^0_{\varphi,T}+(\circledast)^1_{\varphi,T}$.
\end{proof}

\begin{proposition}
\label{scatissweet}
$\cT^\rsc$ is an sw--closed family and consequently $\bqsc$ is iterably
sweet. 
\end{proposition}

\begin{proof}
Plainly $\cT^\rsc$ satisfies the conditions (1) and (2) of \ref{closefam}. 

To verify \ref{closefam}(3) suppose that $\langle T_n:n\leq\omega\rangle
\subseteq\cT^\rsc$ is a sequence of scattered trees such that $\big(\forall
n<\omega\big)\big(T_\omega\cap 2^{\textstyle {\leq}n}= T_n\cap 2^{\textstyle
{\leq} n}\big)$. Let $T=\bigcup\limits_{n\leq\omega} T_n$. We are going to
show that $T$ is a scattered tree, and for this we have to show that $[T]$
is countable. 

Note that if $n<\omega$, $\nu\in\fs\setminus T_\omega$ and $\lh(\nu)\leq n$,
then $\nu\notin T_n$. Therefore, if $\nu\in\fs\setminus T_\omega$ then
$[\nu]\cap [T]\subseteq \bigcup\{[T_n]:n<\lh(\nu)\}$, so $[\nu]\cap [T]$ is
countable. Hence $[T]\setminus [T_\omega]$ is countable and thus (since
$[T_\omega]$ is countable) so is $[T]$.  

The ``consequently'' part follows from \ref{swimsw} (remember that members
of $\cT^\rsc$ are subtrees of $\fs$ so finitely branching).
\end{proof}

Recall that a forcing notion $\bbP$ has $\aleph_1$--caliber if for every
uncountable family $\cF\subseteq\bbP$ there is a condition $p\in\bbP$ such
that $\big|\big\{q\in \cF:q\leq p\big\}\big|=\aleph_1$ (see Truss
\cite{Tr77}). 

\begin{proposition}
\label{nogen}
\begin{enumerate}
\item If a forcing notion $\bbP$ has $\aleph_1$--caliber, then in $\bV^\bbP$
there is no tree $T\subseteq\fs$ such that  
\begin{enumerate}
\item[(a)] for every $\alpha<\omega_1$ there is a countable closed set
  $A\subseteq\can$ coded in $\bV$ such that $\rk(A)=\alpha$ and $A\subseteq
  [T]$, and
\item[(b)] $T$ includes no perfect subtree from $\bV$.
\end{enumerate}
Consequently, $\bbP$ does not add generic object for $\bqsc$.
\item If $\gb>\aleph_1$, then neither the Hechler forcing notion $\bbD$ nor
  its composition $\bbD*\bbC$ with the Cohen real forcing add generic
  objects for $\bqsc$.  
\end{enumerate}
\end{proposition}

\begin{proof}
(1)\quad Suppose toward contradiction that $\bbP$ has an
$\aleph_1$--caliber, $p\in\bbP$ and $\name{T}$ is a $\bbP$--name for a
subtree of $\fs$ such that the condition $p$ forces that both (a) and 
(b) of \ref{nogen}(1) hold true for $\name{T}$. Then for each
$\alpha<\omega_1$ we may choose a scattered tree $T_\alpha\subseteq\fs$ and
a condition $p_\alpha\in\bbP$ such that 
($T_\alpha\in \bV$ and) 
\[\rk([T_\alpha])=\alpha\quad\mbox{ and }\quad p\leq p_\alpha\quad\mbox{ and
}\quad p_\alpha\forces_\bbP\mbox{`` }T_\alpha\subseteq\name{T}\mbox{ ''.}\] 
Since $\bbP$ has an $\aleph_1$--caliber we find a condition $p^*\in\bbP$
such that the set  
\[Y\stackrel{\rm def}{=}\{\alpha<\omega_1:p_\alpha\leq p^*\}\]
is uncountable. Put $T^*=\bigcup\limits_{\alpha\in Y}T_\alpha$. Clearly
$T^*$ is a non-scattered tree and ($T^*\in\bV$ and) $p^*\forces
T^*\subseteq\name{T}$, contradicting (b). 

Concerning the ``consequently'' part it is enough to note that if
$\name{T}^{\rm sc}$ is the canonical $\bqsc$--name for a subset of $\fs$
such that  
\[\forces_{\bqsc}\mbox{`` }\name{T}^{\rm sc}=\bigcup\big\{T^p:p\in
\Gamma_{\bqsc}\big\}\mbox{ ''},\]
then $\forces$`` $\name{T}^{\rm sc}$ is a tree satisfying
\ref{nogen}(1)(a,b) ''.
\medskip

\noindent (2)\quad If the unbounded number $\gb$ is greater than $\aleph_1$,
then both $\bbD$ and $\bbD*\bbC$ have the $\aleph_1$--caliber, so part (1)
applies.  
\end{proof}

\begin{remark}
The forcing notion $\bqsc$ is somewhat similar to the universal forcing
notions discussed in \cite[\S 2.3]{RoSh:672} and \cite{RoSh:845}. However it
follows from \ref{nogen}(2) that if MA holds true, then the composition
$\bbC*\bbD*\bbC$ does not add generic real for $\bqsc$. This is somewhat
opposite to the result presented in \cite[Theorem 2.1]{RoSh:845} and it may
indicate that the answer to the following question is negative.  
\end{remark}

\begin{problem}
Can a finite composition (or, in general, an FS iteration) of the Hechler
forcing notions add a generic object for $\bqsc$ ?
\end{problem}

\section{More sweet examples}
In this section we will present two classes of sw--closed families of trees,
producing many new examples of sweet forcing notions. Let us start with
extending the framework of universality parameters to that of sw--closed
families. 
 
The sweet forcing notions determined by the universality parameters were
introduced in \cite[\S 2.3]{RoSh:672}. In \cite{RoSh:845} we showed that,
unfortunately, the use of them may be somewhat limited because the
composition of, say, the Universal Meager forcing notions adds generic reals
for many examples of the forcing notions determined by universality
parameters. However, as we will show here, {\em families\/} of universality
parameters may determine forcing notions which cannot be embedded into the
known examples of sweet forcing notions.  

Let us start with recalling definitions concerning universality parameters
and the related forcing notions. We will cut down the generality of
\cite[\S 2.3]{RoSh:672} and we will quote here the somewhat simpler setting
of \cite{RoSh:845}. Let $\bH$ be a function from $\omega$ to
$\omega\setminus 2$. 

\begin{definition}
\label{Htree}
\begin{enumerate}
\item {\em A finite $\bH$--tree} is a tree $S\subseteq \bigcup\limits_{n\leq
N}\prod\limits_{i<n}\bH(i)$ with $N<\omega$, $\mrot(S)=\langle\rangle$ and
$\max(S)\subseteq\prod\limits_{i<N}\bH(i)$. The integer $N$ may be called
{\em the level of the tree $S$} and it will be denoted by $\lev(S)$.
\item {\em An infinite $\bH$--tree} is a normal tree $T\subseteq
\bigcup\limits_{n<\omega}\prod\limits_{i<n}\bH(i)$.
\end{enumerate}
\end{definition}

\begin{definition}
\label{univpar}
{\em A simplified universality parameter $\gp$ for $\bH$} is a pair
$(\cG^\gp,F^\gp)=(\cG,F)$ such that   
\begin{enumerate}
\item[$(\alpha)$] elements of $\cG$ are triples $(S,n_\dn,n_\up)$ such that 
$S$ is a finite $\bH$--tree and $n_\dn\leq n_\up\leq\lev(S)$, $(\{\langle
\rangle\},0,0)\in\cG$;  
\item[$(\beta)$] {\bf if:}\qquad $(S^0,n^0_\dn,n^0_\up)\in\cG$, $S^1$ is a
finite $\bH$--tree, $\lev(S^0)\leq\lev(S^1)$, and $S^1\cap 
\prod\limits_{i<\lev(S^0)}\bH(i)\subseteq S^0$, and $n^1_\dn\leq n^0_\dn$,
$n^0_\up\leq n^1_\up\leq\lev(S^1)$,\\  
{\bf then:}\quad $(S^1,n^1_\dn,n^1_\up)\in\cG$, 
\item[$(\gamma)$] $F\in\baire$ is increasing, 
\item[$(\delta)$] {\bf if:}
\begin{itemize}
\item $(S^\ell,n^\ell_\dn,n^\ell_\up)\in\cG$ (for $\ell<2$),
  $\lev(S^0)=\lev(S^1)$,   
\item $S$ is a finite $\bH$--tree, $\lev(S)<\lev(S^\ell)$, and  $S^\ell\cap  
\prod\limits_{i<\lev(S)}\bH(i)\subseteq S$ (for $\ell<2$),   
\item $\lev(S)<n^0_\dn$, $n^0_\up<n^1_\dn$, $F(n^1_\up)<\lev(S^1)$,
\end{itemize}
{\bf then:}\quad there is $(S^*,n^*_\dn,n^*_\up)\in\cG$ such that
\begin{itemize} 
\item $n^*_\dn=n^0_\dn$, $n^*_\up=F(n^1_\up)$, $\lev(S^*)=\lev(S^0)=
\lev(S^1)$, and 
\item $S^0\cup S^1\subseteq S^*$ and $S^*\cap\prod\limits_{i<\lev(S)}\bH(i)=S$. 
\end{itemize}
\end{enumerate}
\end{definition}

\begin{definition}
\label{unforc}
Let $\gp=(\cG,F)$ be a simplified universality parameter for $\bH$. We say
that an infinite $\bH$--tree $T$ is {\em $\gp$--narrow} if  for infinitely
many $n<\omega$, for some $n=n_\dn<n_\up$ we have 
\[(T\cap\bigcup\limits_{n\leq n_\up+1}\prod\limits_{i<n}\bH(i),
n_\dn,n_\up)\in\cG.\]
The family of all $\gp$--narrow infinite $\bH$--trees will denoted by
$\cT^*(\gp,\bH)$.  
\end{definition}

\begin{proposition}
\label{unsw}
If $\gp$ is a simplified universality parameter, then $\cT^*(\gp,\bH)$ is an 
sw--closed family (of finitely branching normal trees). Consequently,
$\bbQ^{\cT(\gp,\bH)}$ is an iterably sweet forcing notion.
\end{proposition}

\begin{proof}
It is should be clear that  $\cT^*(\gp,\bH)$ satisfies
\ref{closefam}(1,2). The proof of \ref{closefam}(3) is, basically, included
in the proof of \cite[Proposition 4.2.5(3)]{RoSh:672}.
\end{proof}

The examples of simplified universality parameters include the following. 

\begin{definition}
[Compare {\cite[Definition 1.7, Example 1.9(2)]{RoSh:845}}]
\label{PPex}
Suppose that the function $\bH$ is increasing and $g\in\baire$ is such that
$(\forall i\in\omega)(0<g(i)<\bH(i))$. Let $A\in\iso$.  We define
$\cG_\bH^{g,A}$ as the family consisting of $(\{\langle\rangle\},0,0)$ and
of all triples $(S,n_\dn,n_\up)$ such that   
\begin{enumerate}
\item[$(\alpha)$] $S$ is a finite $\bH$--tree, $n_\dn\leq n_\up\leq\lev(S)$,
  $A\cap [n_\dn,n_\up]\neq\emptyset$, and 
\item[$(\beta)$] for some sequence $\langle w_i:i\in A\cap[n_\dn, n_\up]
\rangle$ such that $w_i\in [\bH(i)]^{\textstyle {\leq}g(i)}$ (for  $i\in
A\cap[n_\dn, n_\up]$) we have 
\[\big(\forall\eta\in \max(S)\big)\big(\exists i\in A\cap[n_\dn,n_\up)\big)
\big(\eta(i)\in w_i\big).\]    
\end{enumerate}
\end{definition}

\begin{proposition}
\label{PPpart}
Assume that $\bH,g,A$ are as in \ref{PPex}, and $F(n)=\prod\limits_{i\leq
  n}\bH(i)^2$ (for $n\in\omega$). Then $\gp^{g,A}_\bH\stackrel{\rm def}{=}
(\cG^{g,A}_{\bH},F)$ is s simplified universality parameter (and even it is
  a regular universality parameter in the sense of \cite[Definition
  1.14]{RoSh:845}).   
\end{proposition}

The universality parameters $\gp^{g,A}_\bH$ from \ref{PPpart} are related to 
the strong PP--property (see \cite[Ch VI, 2.12*]{Sh:f}, compare also with
\cite[\S 7.2]{RoSh:470}). Note that an infinite $\bH$--tree $T$ is
$\gp^{g,A}_\bH$--narrow if and only if there exist sequences $\bar{w}=
\langle w_i:i\in A\rangle$ and $\bar{n}=\langle n_k:k<\omega\rangle$ such
that    
\begin{itemize}
\item $\big(\forall i\in A\big)\big(w_i\subseteq\bH(i)\ \&\ |w_i|\leq g(i)
\big)$, and 
\item $n_k<n_{k+1}<\omega$ for each $k<\omega$, and 
\item $\big(\forall\eta\in [T]\big)\big(\forall k<\omega\big)\big(\exists
  i\in A\cap [n_k,n_{k+1})\big)\big(\eta(i)\in w_i\big)$.
\end{itemize}
\medskip

It should be clear that the intersection of a family of sw--closed sets of
normal trees is sw--closed. So now we are going to look at the intersections
of the families of $\gp^{g,A}_{\bH}$--narrow trees.

\begin{definition}
Let $\bH,g$ be as in \ref{PPex} and let $\emptyset\neq\cB\subseteq\iso$.
\begin{enumerate}
\item Put $\cT(\cB)=\cT^g_\bH(\cB)\stackrel{\rm def}{=}\bigcap\big\{
\cT^*(\gp^{g,B}_\bH,\bH):B\in \cB\}$ and $\bbP_\cB=\bbQ^{\cT(\cB)}$. 
\item Let $\name{T}_\cB$ be a $\bbP_\cB$--name such that 
\[\forces_{\bbP_\cB}\mbox{`` }\name{T}_\cB=\bigcup\big\{T^p:p\in
\Gamma_{\bbP_\cB}\big\}\mbox{ ''}.\]
\item For a set $A\in\iso$ put 
 \[S_A=\big\{\eta\in\bigcup\limits_{n<\omega}\prod\limits_{i\leq n}\bH(i):
 \big(\forall i\in\lh(\eta)\cap A\big)\big(\eta(i)=0\big)\big\}.\]
\end{enumerate}
\end{definition}

\begin{lemma}
\label{famlem}
Suppose that $\bH,g$ are as in \ref{PPex}.
\begin{enumerate}
\item Let $A,C\in\iso$. Then the tree $S_A$ is $\gp^{g,C}_\bH$--narrow if
  and only if $A\cap C$ is infinite. 
\item Let $\emptyset\neq\cB\subseteq\iso$. Then, in $\bV^{\bbP_\cB}$,
$\name{T}_{\cB}$ is an infinite $\bH$--tree such that   
\begin{enumerate}
\item if $T\in\bV$ is an infinite $\bH$--tree which is
  $\gp^{g,B}_\bH$--narrow for all $B\in\cB$, then there is an $n<\omega$
  such that  
\[\big(\forall\nu\in \name{T}_{\cB}\big)\big(\forall\eta\in T\big)\big(n=
\lh(\nu)<\lh(\eta)\ \Rightarrow\ \nu\conc\eta\rest [n,\lh(\eta))\in
  \name{T}_{\cB}\big),\]  
\item if an infinite $\bH$--tree $T\in\bV$ is not $\gp^{g,B}_\bH$--narrow
  for some $B\in\cB$, then 
\[\big(\forall n<\omega\big)\big(\exists\eta\in T\big)\big(\lh(\eta)>n\ \&\ 
(\forall\nu\in \prod_{i<n}\bH(i))(\nu\conc\eta\rest [n,\lh(\eta))\notin
 \name{T}_{\cB})\big).\]
\end{enumerate}
\end{enumerate}
\end{lemma}

\begin{theorem}
\label{secondnew}
Suppose that $\bbP$ is a ccc forcing notion, $\forces_\bbP$``
$2^{\aleph_0}=\kappa$ '', $\kappa<2^{2^{\aleph_0}}$. Then there is a family
$\cB\subseteq\iso$ such that $\bbP$ does not add the generic object for
the (iterably sweet) forcing notion $\bbP_{\cB}$.
\end{theorem}

\begin{proof}
Note that if $\cU$ is a uniform ultrafilter on $\omega$, then $(\forall
A,B\in\cU)(|A\cap B|=\omega)$ and hence, by \ref{famlem}(1), for every
$A\in\cU$ and every $B\in\cU$, the tree $S_A$ is $\gp^{g,B}_\bH$--narrow.   
Also by \ref{famlem}(1), for every $A\in\iso$ the tree $S_A$ is not 
$\gp^{g,\omega\setminus A}_\bH$--narrow.    

Now, if $\cU',\cU''\subseteq\iso$ are two distinct uniform ultrafilters on 
$\omega$, then we may pick $A\in\iso$ such that $A\in\cU'$ and
$\omega\setminus A\in\cU''$. Then the tree $S_A$
\begin{itemize}
\item is $\gp^{g,B}_\bH$--narrow for every $B\in\cU'$, but   
\item is not $\gp^{g,\omega\setminus A}_\bH$--narrow, $\omega\setminus A\in
  \cU''$. 
\end{itemize}
Therefore, by \ref{famlem}(2), the interpretations of the names
$\name{T}_{\cU'}$, $\name{T}_{\cU''}$ by the corresponding generic filters
must be different. Since there are $2^{2^{\aleph_0}}$ ultrafilters on
$\omega$ we easily get the conclusion.
\end{proof}

\begin{corollary}
There exists an iterably sweet forcing notion $\bbQ$ which cannot be
embedded into the forcing notion constructed in \cite[\S 7]{Sh:176}.
\end{corollary}

Let us present now a different class of sw--closed families of normal trees
and corresponding forcing forcing notions.

\begin{definition}
The {\em sw--closure\/} $\cls(\cT)$ of the family $\cT$ is the smallest
family $\cT^*$ of subtrees of $\fseo$ which includes $\cT$ and is
sw--closed.  
\end{definition}

Clearly, $\cls(\cT)$ is well defined for any family $\cT$ of normal subtrees
of $\fseo$.  

\begin{lemma}
\label{clolem}
\begin{enumerate}
\item Suppose that $T^*$ is a normal subtree of $\fseo$ and let $\cT^*$ be
the family of all normal subtrees of $T^*$. Then $\cT^*$ is sw--closed. 
Consequently, if $\cT\subseteq\cT^*$, then $\cls(\cT)\subseteq\cT^*$.  
\item Assume that $\cT$ is an sw--closed family of normal subtrees of
$\fseo$ and $A\subseteq\baire$ is a closed set. Let 
\[\cT^-(A)=\big\{T\in\cT: [T]\cap A \mbox{ is nowhere dense in }A\big\}.\]
Then $\cT^-(A)$ is sw--closed.
\item If $\cT$ is a family of normal subtrees of $\fseo$, $T\subseteq\fseo$ 
  is a normal tree and 
\[\big(\forall T'\in\cT\big)\big([T]\cap [T']\mbox{ is nowhere dense in }
  [T]\big),\]  
then $T\notin\cls(\cT)$.
\end{enumerate}
\end{lemma}

\begin{proof}
(1)\quad Should be clear.
\medskip

\noindent (2)\quad Clearly $\cT^-(A)$ is closed under finite unions. Assume
now that $T_n,T_\omega\in\cT^-(A)$ are such that $\big(\forall n<\omega\big)
\big(T_\omega\cap \omega^{\textstyle {\leq}n}=T_n\cap\omega^{\textstyle
{\leq} n}\big)$ and let $T=\bigcup\limits_{n\leq\omega} T_n$. We want to
show that $T\in \cT^-(A)$. Since $\cT$ is sw--closed we see that $T\in\cT$, 
so we need to show that $[T]\cap A$ is nowhere dense in $A$. To this end let
$S\subseteq\fseo$ be a normal tree such that $A=[S]$ and suppose that
$\nu\in S$. Since $T_\omega\in \cT^-(A)$, we may find $\eta_0\in S$ such
that $\nu\vartriangleleft\eta_0$ and $\eta_0\notin T_\omega$. Then, by our
assumptions on $\langle T_n:n\leq\omega\rangle$, also for each
$k\geq\lh(\eta_0)$ we have $\eta_0\notin T_k$. Since $T_n\in\cT^-(A)$ (for
$n<\lh(\eta_0)$), the set $\bigcup\limits_{n<\lh(\eta_0)}[T_n]\cap A$ is
nowhere dense in $A$ and hence we may find $\eta\in S$ such that
$\eta_0\vartriangleleft\eta$ and $\eta\notin\bigcup\limits_{n<\lh(\eta_0)}
T_n$. Then we also have $\nu\vartriangleleft\eta\in S$ and $\eta\notin T$. 
\medskip

\noindent (3)\quad Follows from (2).
\end{proof}

\begin{definition}
\begin{enumerate}
\item For a set $A\in\iso$ let $\cT^A$ be the collection of all normal
subtrees $T$ of $\fs$ such that $\big(\forall\nu\in\spliting(T)\big)
\big(\lh(\nu)\in A\big)$.
\item For a family $\cA\subseteq\iso$ let $\cT_\cA=\cls\big(\bigcup
\{\cT^A:A\in\cA\}\big)$.
\end{enumerate}
\end{definition}

\begin{theorem}
\label{newthm}
Suppose that $\bbP$ is a ccc forcing notion, $\forces_\bbP$`` $2^{\aleph_0}
=\kappa$ '', $\kappa<2^{2^{\aleph_0}}$. Then there is a family $\cA
\subseteq\iso$ such that $\bbP$ does not add the generic object for the
(iterably sweet) forcing notion $\bbQ^{\cT_\cA}$. 
\end{theorem}

\begin{proof}
Let us start with some observations of a more general character. 

\begin{claim}
\label{cl4}
Assume that $\cT$ is an sw--closed family of subtrees of $\fs$ such that 
\begin{enumerate}
\item[$(\circledast)$] for every $\eta\in\can$ and $T\in\cT$ we have 
\[T-_2\eta\stackrel{\rm def}{=}\{\nu\in\fs:\nu+_2(\eta\rest\lh(\nu))\in T\}
  \in\cT.\]
\end{enumerate}
Let $\name{T}^{\cT}$ be a $\bqT$--name such that 
\[\forces_{\bqT}\mbox{`` }\name{T}^{\cT}=\bigcup\big\{T^p:p\in\Gamma_{\bqT}
\big\}\mbox{ ''}.\]
Then, in $\bV^{\bqT}$, $\name{T}^{\cT}$ is a subtree of $\fs$ such that 
\begin{enumerate}
\item for every $T\in\cT$ there is an $n<\omega$ such that 

{\em if\/} $\nu_0\in T\cap 2^{\textstyle n}$, $\nu_1\in \name{T}^{\cT}\cap
  2^{\textstyle n}$, and $\nu_0\vartriangleleft\eta\in T$, {\em then\/}
  $\nu_1\conc\eta\rest [n,\lh(\eta))\in \name{T}^{\cT}$, 
\item for every normal tree $T\subseteq\fs$ such that $T\notin\cT$,
$T\in\bV$, we have 
\[\big(\forall n<\omega\big)\big(\exists\eta\in T\big)\big(\lh(\eta)>n\ \&\ 
(\forall\nu\in 2^{\textstyle n})(\nu\conc\eta\rest [n,\lh(\eta))\notin
  \name{T}^{\cT})\big).\]
\end{enumerate}
\end{claim}

\begin{proof}[Proof of the Claim] (1)\quad Suppose that $p\in\bqT$ and
$T\in\cT$. Let $\langle\eta_\ell:\ell<2^{N^p}\rangle$ list all elements of
$\can$ which are constantly zero on $[N^p,\omega)$. It follows from our
assumption $(\circledast)$ that 
\[\big(\forall \ell<2^{N^p}\big)\big(T-_2\eta_\ell\in\cT\ \&\
  T^p-_2\eta_\ell\in\cT\big).\]
Since $\cT$ is sw--closed we may now conclude that (by \ref{closefam}(2)) 
\[T_0\stackrel{\rm def}{=}\bigcup_{\ell<2^{N^p}} (T-_2\eta_\ell)\cup
  \bigcup_{\ell<2^{N^p}} (T^p-_2\eta_\ell)\in\cT,\]
and hence also (by \ref{closefam}(1)) 
\[T_1\stackrel{\rm def}{=}\{\eta\in T_0:(\lh(\eta)\leq N^p\ \&\ \eta\in
T^p)\ \vee\ (\lh(\eta)>N^p\ \&\ \eta\rest N^p\in T^p)\}\in\cT.\]
Now, letting $N^q=N^p$ and $T^q=T_1$ we get a condition $q\in\bqT$ stronger
than $p$ and such that 
\[q\forces(\forall\nu_0\in T\cap 2^{\textstyle N^q})(\forall\nu_1\in
\name{T}^{\cT}\cap 2^{\textstyle N^q})(\forall\eta\in T^{[\nu_0]})(\nu_1
\conc\eta\rest [N^q,\lh(\eta))\in \name{T}^{\cT}).\]
\medskip

\noindent (2)\quad Now suppose that $p\in\bqT$, $n<\omega$ and $T\subseteq
\fs$ is a normal tree which does not belong to $\cT$. Let $N=N^p+n$ and let
$\langle\eta_\ell:\ell<2^N\rangle$ list all elements of $\can$ which are
constantly zero on $[N,\omega)$. It follows from $(\circledast)$ that $T_0 
\stackrel{\rm def}{=}\bigcup\limits_{\ell<2^N} (T^p-_2\eta_\ell)\in\cT$ and
since $T\notin\cT$ we may conclude by \ref{closefam}(1) that $T\setminus
T_0\neq\emptyset$. Pick $\eta\in T\setminus T_0\neq\emptyset$ and note that
necessarily $\lh(\eta)>N\geq n$. Letting $N^q=\lh(\eta)$ and $T^q=T_0$ we
get a condition $q\in\bqT$ stronger than $p$ and such that 
\[q\forces (\forall\nu\in 2^{\textstyle n})(\nu\conc\eta\rest [n,\lh(\eta))
\notin\name{T}^{\cT})\big).\]
\end{proof}

\begin{claim}
\label{cl5}
If $\cT$ is a collection of normal subtrees of $\fs$ such that the demand in
\ref{cl4}$(\circledast)$ holds for $\cT$, then also $\cls(\cT)$ satisfies
this condition. Consequently, for each $\cA\subseteq\iso$, $(\circledast)$
of \ref{cl4} holds true for $\cT_{\cA}$.
\end{claim}

\begin{proof}[Proof of the Claim]
Should be clear.
\end{proof}

\begin{claim}
\label{cl3}
Suppose that $A\in\iso$ and $\cA\subseteq\iso$ are such that 
\[\big(\forall B\in\cA\big)\big(|A\setminus B|=\omega).\]
Then $\cT^A\nsubseteq\cT_{\cA}$.
\end{claim}

\begin{proof}[Proof of the Claim]
Let $T=\{\nu\in\fs:(\forall n<\lh(\nu))(\nu(n)=1\ \Rightarrow\ n\in
A)\}$. Plainly $T\in \cT^A$. Also, for every $B\in\cA$ and $T'\in\cT^B$ the
set $[T]\cap [T']$ is nowhere dense in $[T]$, so by \ref{clolem}(3) $T\notin
\cls\big(\bigcup\{\cT^B:B\in\cA\}\big)=\cT_{\cA}$.
\end{proof}

Now choose a family $\cI\subseteq\iso$ of almost disjoint sets, $|\cI|=
2^{\aleph_0}$. 

Suppose that $\cA,\cB\subseteq\cI$, $\cA\neq\cB$, say $A\in\cA\setminus
\cB$. Then $(\forall B\in\cB)(|A\setminus B|=\omega)$ and hence (by Claim
\ref{cl3}) we get $\cT^A\nsubseteq\cT_{\cB}$, so we have a normal tree $T\in 
\cT_\cA\setminus\cT_\cB$. Now look at Claim \ref{cl4} --- by \ref{cl5} it is
applicable to $\bbQ^{\cT_\cA}$, $\bbQ^{\cT_\cB}$ and we get from it that if
$T_\cA,T_\cB\subseteq\fs$ are trees generic over $\bV$ for $\bbQ^{\cT_\cA}$,
$\bbQ^{\cT_\cB}$, respectively, then  
\begin{itemize}
\item $\big(\exists n<\omega\big)\big(\forall\nu\in T_\cA\cap 2^{\textstyle
  n}\big)\big(\forall\eta\in T\big)\big(\lh(\eta)>n\ \Rightarrow\ \nu\conc\eta
  \rest [n,\lh(\eta))\in T_\cA\big)$,  
\item $\big(\forall n<\omega\big)\big(\exists\eta\in
  T\big)\big(\lh(\eta)>n\ \&\  (\forall\nu\in 2^{\textstyle
  n})(\nu\conc\eta\rest [n,\lh(\eta))\notin T_\cB\big)$.
\end{itemize}
Hence $T_\cA\neq T_\cB$. Since $\bbP$ satisfies the ccc and $\forces_\bbP$``
$2^{\aleph_0}=\kappa$ '' and $\kappa<2^{2^{\aleph_0}}$, we may find a family
$\cF$ of subsets of $\cI$ such that $|\cF|=\kappa$ and 
\[\forces_\bbP
\mbox{`` for no $\cA\subseteq\cI$ with $\cA\notin\cF$, there is a
  $\bbQ^{\cT_\cA}$--generic filter over $\bV$ ''.}\]
\end{proof}

One should note that the examples of sweet forcing notions which cannot be
embedded into the one constructed in \cite[\S 7]{Sh:176} which we gave in
this section are not very nice --- it may well be that the parameters
$\cA,\cB$ needed to define them are not definable from a real. Even the
candidate for a somewhat definable example from the previous section, the
forcing notion $\bqsc$, is not Souslin. Thus the following variant of
\cite[Problem 5.5]{RoSh:672} may be of interest. 

\begin{problem}
Is there a Souslin ccc iterably sweet forcing notion $\bbQ$ such that no
finite composition of the Universal Meager forcing notion adds a
$\bbQ$--generic real? Such that the forcing of \cite[\S 7]{Sh:176} does not
add $\bbQ$--generic real?
\end{problem}
 
\section{Subforcings, Quotients and likes}
Topological sweetness, as defined in \ref{topsweet}, is a property of
particular representation of a forcing notion. It is only natural to ask
if a forcing notion having a topologically sweet dense subforcing is
topologically sweet, or, in general, if a forcing notion equivalent to a
topologically sweet one is topologically sweet. We start this section with
some results in these directions.

\begin{definition}
\label{glb}
We say that a forcing notion $\bbP$ has a GLB--property provided that for
every $p_0,\ldots,p_k\in\bbP$, $k<\omega$, there is $q\in\bbP$ such that 
\begin{enumerate}
\item[$(\alpha)$] $q\leq p_i$ for $i\leq k$, and 
\item[$(\beta)$]  if $q^*\in\bbP$ satisfies $(\forall i\leq k)(q^*\leq
  p_i)$, then $q^*\leq q$.
\end{enumerate}
\end{definition}

\begin{remark}
If $\BB$ is a Boolean algebra, then $\BB^+$ is a forcing notion with the
GLB--property. Also the forcing notions $\mbR$ and $\bbA$ defined in
\ref{amoran} later have this property. 
\end{remark}

\begin{proposition}
\label{transfer}
Suppose that a forcing notion $\bbP$ has the GLB--property and $\bbQ
\subseteq\bbP$ is its dense subforcing. If $\bbQ$ is topologically sweet,
then so is $\bbP$.
\end{proposition}

\begin{proof}
Let $(\bbQ,\cB)$ be a model of topological sweetness and let $\tau$ be the
topology on $\bbQ$ generated by $\cB$. For sets $U_0,\ldots,U_k\in\cB$,
$k<\omega$, define
\[W(U_0,\ldots,U_k)=\{p\in\bbP:(\forall i\leq k)(\exists q\in U_i)(p\leq
q)\},\] 
and let 
\[\cB^*=\big\{W(U_0,\ldots,U_k):k<\omega\ \&\ U_0,\ldots,U_k\in\cB\big\}
\cup\big\{\{\bO_\bbP\}\big\}.\]
It should be clear that 
\begin{itemize}
\item $\cB^*$ is closed under finite intersections, and
\item it is a countable basis of a topology $\tau^*$ on $\bbP$, and
\item $\bO_\bbP$ is an isolated point in $\tau^*$.
\end{itemize}
We are going to show that the topology $\tau^*$ satisfies the demand of
\ref{topsweet}(ii). So suppose that a sequence $\bar{p}=\langle p_n:n<
\omega\rangle\subseteq\bbP$ is $\tau^*$--converging to $p\in\bbP$ and $q\geq
p$ and $W$ is a $\tau^*$--neighbourhood of $q$. Pick $U_0,\ldots, U_k\in\cB$
such that $q\in W(U_0,\ldots, U_k)\subseteq W$ and let $q_i\in U_i$ (for
$i\leq k$) be such that $q\leq q_i$. Furthermore, for $i\leq k$, let
$\{V^i_n:n<\omega\}$ be a basis of $\tau$--neighbourhoods of $q_i\in\bbQ$
such that $(\forall n_0<n_1<\omega)(q_i\in V^i_{n_1}\subseteq V^i_{n_0}
\subseteq U_i)$.  

Since $p\in W(V^0_n,V^1_n,\ldots,V^k_n)\in\cB^*$ (for each $n<\omega$) and
the sequence $\bar{p}$ $\tau^*$--converges to $p$, we may choose an
increasing sequence $\langle m_n:n<\omega\rangle\subseteq\omega$ such that
$\big(\forall n<\omega\big)\big(p_{m_n}\in W(V^0_n,V^1_n,\ldots,V^k_n)
\big)$. Then we may also pick $p^*_{n,i}$ (for $n<\omega$ and $i\leq k$)
such that $p_{m_n}\leq p^*_{n,i}\in V^i_n$. Fix $i\leq k$ and look at the
sequence $\bar{p}^*_i=\langle p^*_{n,i}:n<\omega\rangle$: clearly it
$\tau$--converges to $q_i$. Consequently, we may easily choose (be repeated
application of \ref{topsweet}(ii) for $\tau$) conditions $q^*_i\in\bbQ$ such 
that 
\begin{itemize}
\item $q_i\leq q^*_i\in U_i$ for $i\leq k$, and 
\item $(\exists^\infty n<\omega)(\forall i\leq k)(p^*_{n,i}\leq q^*_i)$. 
\end{itemize}
Since $\bbP$ has the GLB--property we may pick $q^*\in\bbP$ such that 
\begin{enumerate}
\item[$(\alpha)$] $q^*\leq q_i$ for $i\leq k$, and 
\item[$(\beta)$]  if $r\in\bbP$ is weaker than $q^*_0,\ldots,q^*_k$, then
  $r\leq q^*$. 
\end{enumerate}
Then, plainly, $q^*\in W(U_0,\ldots,U_k)$ and $q\leq q^*$ and
$(\exists^\infty n<\omega)(p_{m_n}\leq q)$. 
\end{proof}

\begin{proposition}
\label{closedown}
Assume that $\bbP$ is a topologically sweet forcing notion. Then there is a
model $(\bbP,\cB^*)$ of topological sweetness such that all members of
$\cB^*$ are downward closed.  
\end{proposition}

\begin{proof}
Let $(\bbP,\cB)$ be a model of topological sweetness. For $U\in \cB$ put
$W(U)=\{p\in\bbP:(\exists q\in U)(p\leq q)\}$, and let
$\cB^*=\{W(U):U\in\cB\}$. Note that if $p\in W(U_0)\cap W(U_1)$ and $p\leq
p_0\in U_0$, $p\leq p_1\in U_1$, then there is $V\in\cB$ such that $p\in V$
and $V\subseteq W(U_0)\cap W(U_1)$ (remember \ref{bastoplem}(1)). Hence we
easily conclude that $\cB^*$ is a base of a topology $\tau^*$ on
$\bbP$. Similarly as in \ref{transfer} one shows that $(\bbP,\cB^*)$ is a
model of topological sweetness.    
\end{proof}

\begin{proposition}
\label{tosub}
Assume that $\bbP$ is a topologically sweet and separative partial order,
$\bbQ$ is a forcing notion. Suppose also that  
\[(\forall q\in\bbQ)(\exists p\in\bbP)(p\forces_\bbP\mbox{`` there is a
  $\bbQ$--generic $\name{G}\subseteq\bbQ$ over $\bV$ such that
  $q\in\name{G}$ ''}).\] 
Then $\bbQ$ is equivalent to a topologically sweet forcing notion.  
\end{proposition}

\begin{proof}
It follows from our assumptions on $\bbP$ that it is (isomorphic to) a dense
subset of ${\bf BA}(\bbP)^+$ and hence, by \ref{transfer}+\ref{closedown},
there is a model $({\bf BA}(\bbP)^+,\cB)$ of topological sweetness such that
all members of $\cB$ are downward closed. By the assumptions on $\bbQ,\bbP$
we also know that ${\bf BA}(\bbQ)$ is a complete subalgebra of ${\bf
  BA}(\bbP)$; let $\pi:{\bf BA}(\bbP)\longrightarrow{\bf BA}(\bbQ)$ be the
projection. Put 
\[\cB'=\{U\cap{\bf BA}(\bbQ)^+:U\in\cB\}.\]
We claim that $({\bf BA}(\bbQ)^+,\cB')$ is a model of topological
sweetness. It is easy to verify \ref{topsweet}(i), so let us only argue that
\ref{topsweet}(ii) holds true. To this end suppose that a sequence $\bar{p}=
\langle p_n:n<\omega\rangle\subseteq{\bf BA}(\bbQ)^+$ converges to $p\in
{\bf BA}(\bbQ)^+$ (in the topology generated by $\cB'$) and let $p\leq q\in
U\cap {\bf BA}(\bbQ)^+$, $U\in\cB$. Then also $\bar{p}$ converges to $p$ in
the topology generated by $\cB$ on ${\bf BA}(\bbP)^+$, so we may find $r\in
{\bf BA}(\bbP)^+$ such that $q\leq r\in U$ and $(\exists^\infty n<\omega)(
p_n\leq r)$. Let $r^*=\pi(r)\in {\bf BA}(\bbQ)$. Then we have
\begin{itemize}
\item $q\leq r^*$ (as $\pi$ is the projection and $q\in {\bf BA}(\bbQ)^+$,
  $q\leq r$),
\item $(\exists^\infty n<\omega)(p_n\leq r^*)$ (as $\pi$ is the projection
  and $p_n\in {\bf BA}(\bbQ)^+$), 
\item $r^*\in U$ (as $U$ is downward closed, $r^*\leq r\in U$).
\end{itemize}
\end{proof}

The sweetness and topological sweetness are important properties because
they are preserved in amalgamations of forcing notions. Since the
amalgamation can be represented as the composition with the product of two
quotients (see, e.g., \cite{JuRo} on that), one may ask if sweetness is also
preserved in quotients. 

\begin{definition}
Let $\bbP,\bbQ$ be forcing notions and suppose that $\bbQ\lesdot
{\bf BA}(\bbP)$. The quotient $(\bbP:\bbQ)$ is the {\em $\bbQ$--name\/} for
the subforcing of $\bbP$ consisting of all $p\in\bbP$ such that $p$ is
compatible (in ${\bf BA}(\bbP)$) with all members of $\Gamma_{\bbQ}$. Thus
for $p\in\bbP$ and $q\in\bbQ$,
\[\begin{array}{l}
q\forces_{\bbQ}\mbox{`` }p\in (\bbP:\bbQ)\mbox{ ''\qquad if and only if
  \qquad }\\
(\forall r\in\bbQ)(q\leq r\ \Rightarrow\ r,p\mbox{ are compatible
  in } {\bf BA}(\bbP)).
  \end{array}\]
\end{definition}

\begin{theorem}
\label{overc}
Let $\bbC$ be the standard Cohen forcing notion (so it is a countable
atomless partial order). Suppose that $(\bbP,\cB)$ is a model of topological
sweetness and $\bbC\lesdot {\bf BA}(\bbP)$. Let $\name{\cB}^\bbC$ be the
$\bbC$-name for the family $\{U\cap (\bbP:\bbC): U\in\cB\}$. Then 
\[\forces_\bbC\mbox{`` }\big((\bbP:\bbC),\name{\cB}^\bbC\big)\mbox{ is
  a model of topological sweetness ''}.\] 
\end{theorem}

\begin{proof}
First note that, in $\bV^{\bbC}$, $\name{\cB}^\bbC$ is a countable basis of
a topology on $(\bbP:\bbC)$, and $\bO_{(\bbP:\bbC)}=\bO_{\bbP}$ is an
isolated point in this topology. Thus the only thing that we should verify
is the demand in \ref{topsweet}(1)(ii). 

Suppose that $\eta\in\bbC$ and $\bbC$--names $\langle \name{p}_i:i<\omega
\rangle$, $\name{p},\name{q}$ and $\name{W}$ are such that 
\[\begin{array}{ll}
\eta\forces_{\bbC}&\mbox{`` }\name{p}_i,\name{p},\name{q}\in (\bbP:\bbC),\
\name{W}\in \name{\cB}^{\bbC},\ \name{p}\leq \name{q}\in \name{W}\mbox{
  and}\\
&\mbox{ the sequence } \langle \name{p}_i:i<\omega\rangle\mbox{ converges to
  $\name{p}$ in the topology generated by $\name{\cB}^\bbC$ ''}
  \end{array}\]
Passing to a stronger than $\eta$ condition in $\bbC$ (if necessary), we may
assume that for some $p,q\in\bbP$ and $W\in \cB$ we have
\[\eta\forces_\bbC\mbox{`` }\name{p}=p\ \&\ \name{q}=q\ \&\ \name{W}=W\cap
(\bbP:\bbC)\mbox{ ''.}\]
Then also $\eta\forces_{\bbC}$`` $p,q\in(\bbP:\bbC)$ '' and $p\leq q\in
W$. Let us choose a condition $q^+\in\bbP$ which is (in ${\bf BA}(\bbP)$)
stronger than both $q$ and $\eta$, and let $U\in\cB$ be a neighborhood of
$q^+$ such that any two members of $U$ are compatible in $\bbP$ (remember
\ref{bastoplem}(2)). Next, choose $W^+\in \cB$ such that $q\in W^+\subseteq
W$ and every member of $W^+$ has an upper bound in $U$ (possible by 
\ref{bastoplem}(1)). 

Pick $V_i\in\cB$ (for $i<\omega$) such that $\{V_i:i<\omega\}$ forms a
neighbourhood basis at $p$ (for the topology generated by $\cB$) such that
for each $i<\omega$: 
\begin{enumerate}
\item[$(\alpha)$] $p\in V_{i+1}\subseteq V_i$,    
\item[$(\beta)$]  any $i+1$ elements of $V_{i+1}$ have a common upper bound
  in  $V_i$. 
\end{enumerate}
[The choice is clearly possible; remember \ref{bastoplem}.] 

Clearly $\eta\forces_{\bbC}\mbox{`` }\{V_i\cap (\bbP:\bbC):i<\omega\}$
forms a neighbourhood basis at $p$ (for the topology generated by
$\name{\cB}^\bbC$) ''. Hence, without loss of generality, we may assume that 
$\eta\forces_{\bbC}$`` $\name{p}_i\in V_i$ '' (as we may change the names
$\name{p}_i$ reflecting a passage to a subsequence). Let us fix a list
$\{\nu_\ell:\ell<\omega\}$ of all conditions in $\bbC$ stronger than $\eta$,
and for every $i,\ell<\omega$ let us pick $p_{i,\ell}\in\bbP$ such that
$\nu_\ell\not\forces_{\bbC}\mbox{`` }\name{p}_i\neq p_{i,\ell}\mbox{
  ''}$. Note that then $p_{i,\ell}\in V_i$, so by clause $(\beta)$ above we
may choose $p^*_i\in V_i$ such that for each $i>0$ we have 
\[(\forall \ell\leq i)(p_{i+1,\ell}\leq p^*_i).\]
The sequence $\langle p_i^*:i<\omega\rangle$ converges to $p$ so (by
\ref{topsweet}(1)(ii) for $(\bbP,\cB)$) there are a condition $r\in\bbP$ and
an infinite set $A\subseteq\omega$ such that 
\[r\in W^+\ \mbox{ and }\ q\leq r\ \mbox{ and }\ (\forall i\in A)(p^*_i \leq
r).\] 
By the choice of $W^+$, the condition $r$ has an upper bound in $U$ and
hence (by the choice of $U$) $r,q^+$ are compatible in $\bbP$. Therefore, as
$q^+$ is stronger than $\eta$ (in ${\bf BA}(\bbP)$), there is $\nu\in\bbC$
stronger than $\eta$ such that $\nu\forces_{\bbC}$`` $r\in(\bbP:\bbC)$ ''.  
Now the proof follows from the following Claim.

\begin{claim}
\label{cl1}
$\nu\forces_{\bbC}$`` $(\exists^\infty i<\omega)(\name{p}_i\leq r)$ ''.
\end{claim}

\begin{proof}[Proof of the Claim]
If not, then we may find $\nu'\in\bbC$ stronger than $\nu$ and $i'<\omega$
such that $\nu'\forces_{\bbC}$`` $(\forall i\geq i')(\name{p}_i\nleq r)$ ''. 
Let $\ell<\omega$ be such that $\nu'=\nu_\ell$ and let $i\in A$ be larger
than $\ell+i'+1$. Look at our choices before - we know that:
\begin{enumerate}
\item[(i)]   $p^*_i\leq r$,
\item[(ii)]  $p_{i+1,\ell}\leq p^*_i$,  
\item[(iii)] $\nu_\ell\not\forces_{\bbC}\mbox{`` }\name{p}_{i+1}\neq p_{i+1, 
\ell}\mbox{ ''}$.
\end{enumerate}
Therefore some condition $\nu^*\in\bbC$ stronger than $\nu_\ell$ forces that
$\name{p}_{i+1}\leq r$, contradicting the choice of $\nu'=\nu_\ell$ (as
$i+1>i'$).  
\end{proof}
\end{proof}

In the rest of this section we are going to show that the result of
\ref{overc} cannot be very much improved: when taking a quotient over a
random real forcing we may loose topological sweetness. Let us start with
recalling some notation and definitions, which we will need later. 

\begin{definition}
\label{amoran}
\begin{enumerate}
\item  The Lebesgue (product) measure on $\can$ is denoted by $\mul$,
$\borel(\can)$ is the $\sigma$--field of Borel subsets of $\can$ and $\bbL$
is the $\sigma$--ideal of Lebesgue null subsets of $\can$. The quotient
complete Boolean algebra $\BB=\borel(\can)/\bbL$ is called {\em the random
algebra}.  
\item The random forcing notion $\mbR$ is defined as follows:\\
{\bf a condition} in $\mbR$ is a closed subset of $\can$ of positive
Lebesgue measure,\\
{\bf the order} of $\mbR$ is the reverse inclusion. 
\item The amoeba for measure forcing notion $\bbA$ is defined as follows:\\
{\bf a condition} in $\bbA$ is a closed subset $F$ of $\can$ such that
$\mul(F)>\frac{1}{2}$,\\
{\bf the order} of $\bbA$ is the reverse inclusion. 
\end{enumerate}
\end{definition}

Of course, $\BB={\bf BA}(\mbR)$. Let us also recall that both $\mbR$ and
$\bbA$ are topologically sweet (see \cite[1.3.3]{St85}).

\begin{proposition}
\label{amoeba}
\begin{enumerate}
\item $\forces_{\BB}$`` $\bbA^\bV$ is not topologically sweet ''.
\item $\forces_{\BB}$`` $\mbR^\bV$ is not topologically sweet ''.
\end{enumerate}
\end{proposition}

\begin{proof}
(1)\quad Suppose toward contradiction that 
\[\lbv\mbox{ there is a model of topological sweetness based on
}\bbA^\bV\rbv_\BB\neq {\bf 0}_\BB.\]
Since the random algebra is homogeneous, we may assume that we have
$\BB$--names $\name{U}_n$ for subsets of $\bbA^\bV$ such that 
\begin{enumerate}
\item[$(*)_0$] $\forces_\BB$`` $(\bbA^\bV,\{\name{U}_n:n<\omega\})$ is a
  model of topological sweetness ''.
\end{enumerate}
For $i<\omega$ let $m_i=\lfloor -\frac{i\cdot 2^i}{\log_2(1-2^{-2^{i+1}})}
\rfloor+2$, so $\frac{m_i}{2^i}>\frac{-i}{\log_2(1-2^{-2^{i+1}})}$ and thus 
\begin{enumerate}
\item[$(*)_1$] $(1-2^{-2^{i+1}})^{m_i/2^i}<2^{-i}$.
\end{enumerate}
Let $\mu$ be the product Lebesgue measure on the space
$\prod\limits_{i<\omega} m_i$, and let $\mu^*$ be the corresponding outer
measure.
  
Define $\langle n_i:i<\omega\rangle$ by $n_0=0$, $n_{i+1}=n_i+m_i\cdot
2^{i+1}$, and for $i<\omega$, $j<m_i$ put 
\[t^i_j\stackrel{\rm def}{=}\big\{\sigma\in 2^{\textstyle [n_i,n_{i+1})}:
\big(\exists\ell<2^{i+1}\big)\big(\sigma(n_i+j\cdot 2^{i+1}+\ell)=1\big)
\big\}.\] 
Note that 
\begin{enumerate}
\item[$(*)^i_2$] if $j_0<j_1<\ldots<j_k<m_i$, then
\[\big|t^i_{j_0}\cap t^i_{j_1}\cap\ldots\cap t^i_{j_k}\big|=\big(1-
2^{-2^{i+1}}\big)^{k+1}\cdot 2^{m_i\cdot 2^{i+1}}.\]
\end{enumerate}
For $x\in\prod\limits_{i<\omega}m_i$ let 
\[Z_x\stackrel{\rm def}{=}\big\{\eta\in\can:\big(\forall i<\omega\big)
\big(\eta\rest [n_i,n_{i+1})\in t^i_{x(i)}\big)\big\}\]
and note that $Z_x$ is a closed set and $\mul(Z_x)>\frac{1}{2}$, so
$Z_x\in\bbA$. For each $x\in\prod\limits_{i<\omega}m_i$ and $n<\omega$ we
may pick a Borel set $B(x,n)\subseteq\can$ such that $\lbv Z_x\in
\name{U}_n\rbv_\BB=[B(x,n)]_\bbL$. Next, for each $k<\omega$ (and
$x\in\prod\limits_{i<\omega}m_i$ and $n<\omega$) choose a clopen set
$C(x,n,k)\subseteq\can$ such that $\mul\big(B(x,n)\vartriangle
C(x,n,k)\big)<2^{-k}$. Now, for $n<\omega$, consider a binary relation
$\sim_n$ on $\prod\limits_{i<\omega}m_i$ given by
\[x\sim_n y\quad\mbox{ if and only if }\quad \big(\forall k,\ell\leq n\big)
\big(C(x,\ell,k)=C(y,\ell,k)\big).\]
It should be clear that (for each $n<\omega$) $\sim_n$ is an equivalence
relation on $\prod\limits_{i<\omega} m_i$ such that 
\begin{enumerate}
\item[$(*)^n_3$] $x\sim_{n+1} y\ \Rightarrow\ x\sim_n y$ \quad (for each
  $x,y\in\prod\limits_{i<\omega} m_i$), and  
\item[$(*)^n_4$] $\prod\limits_{i<\omega}m_i/\sim_n$ is countable.
\end{enumerate}
Consequently we may pick $x^*\in\prod\limits_{i<\omega}m_i$ such that for
each $n<\omega$ we have 
\[\lim_{\ell\to\infty}\frac{\mu^*\big(\big\{x\in\prod\limits_{i<\omega}m_i:
x\rest\ell=x^*\rest\ell\ \ \&\ \ x\sim_n x^*\big\}\big)}{\mu\big(\big\{x\in 
\prod\limits_{i<\omega}m_i:x\rest\ell=x^*\rest\ell\big\}\big)}=1.\]
So now we may choose an increasing sequence $\langle\ell_i:i<\omega\rangle
\subseteq\omega$ such that for $i<\omega$ we have
\[\mu^*\big(\big\{x\in\prod\limits_{j<\omega}m_j:x\rest\ell_i=x^*\rest
\ell_i\ \ \&\ \ x\sim_i x^*\big\}\big)>\frac{1}{2}\mu\big(\big\{x\in
\prod\limits_{j<\omega}m_j:x\rest\ell_i=x^*\rest\ell_i\big\}\big),\]
and then for each $i<\omega$ we may choose $v_i\subseteq m_{\ell_i}$ and
$\langle y^i_k:k\in v_i\rangle\subseteq\prod\limits_{j<\omega} m_j$ such
that 
\begin{enumerate}
\item[$(*)^i_5$] $|v_i|>\frac{1}{2} m_{\ell_i}$, 
\item[$(*)^i_6$] $y^i_k\rest \ell_i=x^*\rest \ell_i$, $y^i_k(\ell_i)=k$ and
  $y^i_k\sim_i x^*$ for $k\in v_i$.
\end{enumerate}
It follows from the definition of the relations $\sim_n$ and from $(*)^i_6$
that for each $k\in v_i$ and all $\ell\leq i$ we have 
\[\mul\big(B(x^*,\ell)\vartriangle B(y^i_k,\ell)\big)<2^{1-i}.\]
Thus, for each $i<\omega$, we may pick a partition $\langle B^i_k:k\in
v_i\rangle$ of $\can$ into disjoint Borel sets such that for all $k\in v_i$
we have 
\begin{enumerate}
\item[$(*)^{i,k}_7$] $\mul(B^i_k)=\frac{1}{|v_i|}$, and 
\item[$(*)^{i,k}_8$] $\mul\big(B^i_k\cap(B(x^^*,\ell)\vartriangle B(y^i_k,
  \ell))\big)<2^{1-i}/|v_i|$ for all $\ell\leq i$.
\end{enumerate}
Let $\name{x}_i$ be a $\BB$--name for a member of $\bV\cap
\prod\limits_{j<\omega} m_j$ such that 
\[\big(\forall k\in v_i\big)\big(\lbv\name{x}_i=y^i_k\rbv_\BB=[B^i_k]_\bbL
\big).\] 

\begin{claim}
\label{cl2}
\[\forces_\BB\mbox{`` }\big(\forall n<\omega\big)\big(\forall^\infty i<\omega
\big)\big(Z_{x^*}\in\name{U}_n\ \Rightarrow\ Z_{\name{x}_i}\in
\name{U}_n\big).\]
\end{claim}

\begin{proof}[Proof of the Claim]
Note that for $n,i<\omega$ we have
\[\lbv Z_{\name{x}_i}\notin \name{U}_n\rbv_\BB=\Big[\bigcup_{k\in v_i} B^i_k 
\setminus B(y^i_k,n)\Big]_\bbL,\]
and thus $\lbv Z_{x^*}\in\name{U}_n\ \&\ Z_{\name{x}_i}\notin\name{U}_n\rbv_\BB
=\big[\bigcup\limits_{k\in v_i}\big(B(x^*,n)\setminus B(y^i_k,n)\big)\cap
B^i_k\big]_\bbL$. It follows from $(*)^{i,k}_8$ that  (for $n\leq i<\omega$)
we have 
\[\mul\big(\bigcup_{k\in v_i}\big(B(x^*,n)\setminus B(y^i_k,n)\big)\cap
B^i_k\big)<2^{1-i}.\]
Hence for, each $n<\omega$,
\[\mul\Big(\bigcap_{m<\omega}\bigcup_{i>m}\big(\bigcup_{k\in v_i}
\big(B(x^*,n)\setminus B(y^i_k,n)\big)\cap B^i_k\big)\Big)=0,\]
so 
\[\lbv (\exists^\infty i<\omega)(Z_{x^*}\in\name{U}_n\ \&\
Z_{\name{x}_i}\notin \name{U}_n\rbv_\BB={\bf 0}_\BB,\]
and the Claim follows.
\end{proof}

It follows from $(*)_0$ and \ref{cl2} that 
\[\lbv\big(\exists F\in\bbA^\bV\big)\big(F\subseteq Z_{x^*}\ \&\
(\exists^\infty i<\omega)(F\subseteq Z_{\name{x}_i})\big)\rbv_\BB={\bf
  1}_\BB,\] 
and therefore we may find $F\in\bbA\cap\bV$ such that $F\subseteq Z_{x^*}$
and $a\stackrel{\rm def}{=}\lbv(\exists^\infty i<\omega)(F\subseteq
Z_{\name{x}_i})\rbv_\BB\neq{\bf 0}_\BB$. For $i<\omega$ put 
\[w_i=\{k\in v_i:F\subseteq Z_{y^i_k}\}\quad \mbox{ and }\quad C_i=
  \bigcup_{k\in w_i} B^i_k.\]
Plainly, $a=\big[\bigcap\limits_{m<\omega}\bigcup\limits_{i>m} C_i
\big]_\bbL$ so (as $a\neq{\bf 0}_\BB$) $\sum\limits_{i=1}^\infty\mul(C_i)=
\infty$, and hence the set 
\[I\stackrel{\rm def}{=}\{i<\omega:\mul(C_i)>2^{1-i}\}\]
is infinite. 

Fix $i\in I$ for a moment. Then 
\[2^{1-i}<\mul(C_i)=\sum_{k\in w_i}\mul(B^i_k)=\frac{|w_i|}{|v_i|},\]
and thus (by $(*)^i_5$)
\[|w_i|>|v_i|\cdot 2^{1-i}>\frac{1}{2}\cdot m_{\ell_i}\cdot 2^{1-i}\geq
m_{\ell_i}/2^{\ell_i}.\] 
Hence, by $(*)^{\ell_i}_1$, we get $\big(1-2^{-2^{\ell_i+1}}\big)^{|w_i|}<
2^{-\ell_i}\leq 2^{-i}$. Now (for our $i\in I$) consider the closed set
$Y_i\stackrel{\rm def}{=}\bigcap\limits_{k\in w_i} Z_{y^i_k}$ and note that 
\[Y_i\subseteq\big\{\eta\in\can: (\forall k\in w_i)(\eta\rest [n_{\ell_i},
  n_{\ell_i+1})\in t^{\ell_i}_k)\big\}.\]
Thus, by $(*)^{\ell_i}_2$, we may conclude that (for our $i\in I$) 
\[\mul(Y_i)\leq\frac{\big|\bigcap\limits_{k\in w_i} t^{\ell_i}_k\big|}{
2^{m_{\ell_i}\cdot 2^{\ell_i+1}}}=\big(1-2^{-2^{\ell_i+1}}\big)^{|w_i|}<
2^{-i}.\]

Since $I$ is infinite and for every $i\in I$ we have $F\subseteq
\bigcap\limits_{k\in w_i} Z_{y^i_k}=Y_i$ we may now conclude that
$\mul(F)=0$, contradicting $F\in\bbA$. 
\medskip

\noindent (2)\quad The same proof as for (1) works here too. 
\end{proof}

Putting together \ref{transfer} and \ref{amoeba} we may easily conclude the
following.  

\begin{corollary}
\label{overran}
Both $\mbR\times\mbR$ and $\bbA\times\bbA$ are topologically sweet, but 
\[\begin{array}{ll}
\forces_\mbR&\mbox{`` no dense subforcing of $(\mbR\times\mbR:\mbR)$
  ($(\mbR\times\bbA:\mbR)$, respectively)}\\
&\quad\mbox{is topologically sweet ''.}
  \end{array}\]

\end{corollary}


\bigskip\bigskip\bigskip

\begin{thebibliography}{10}

\bibitem{BaJu95}
Tomek Bartoszy\'nski and Haim Judah.
\newblock {\em {Set Theory: On the Structure of the Real Line}}.
\newblock A K Peters, Wellesley, Massachusetts, 1995.

\bibitem{J}
Thomas Jech.
\newblock {\em {Set theory}}.
\newblock Academic Press, New York, 1978.

\bibitem{JuRo}
Haim Judah and Andrzej Ros{\l}anowski.
\newblock {On Shelah's Amalgamation}.
\newblock In {\em Set Theory of the Reals}, volume~6 of {\em Israel
  Mathematical Conference Proceedings}, pages 385--414. 1992.

\bibitem{RoSh:845}
Andrzej Roslanowski and Saharon Shelah.
\newblock {Universal forcing notions and ideals}.
\newblock {\em Journal of Symbolic Logic}, submitted.
\newblock math.LO/0404146.

\bibitem{RoSh:470}
Andrzej Roslanowski and Saharon Shelah.
\newblock {Norms on possibilities I: forcing with trees and creatures}.
\newblock {\em {Memoirs of the American Mathematical Society}}, 141(671):xii +
  167, 1999.
\newblock math.LO/9807172.

\bibitem{RoSh:672}
Andrzej Roslanowski and Saharon Shelah.
\newblock {Sweet {\&} Sour and other flavours of ccc forcing notions}.
\newblock {\em Archive for Mathematical Logic}, 43:583--663, 2004.
\newblock math.LO/9909115.

\bibitem{Sh:176}
Saharon Shelah.
\newblock {Can you take Solovay's inaccessible away?}
\newblock {\em {Israel Journal of Mathematics}}, 48:1--47, 1984.

\bibitem{Sh:f}
Saharon Shelah.
\newblock {\em {Proper and improper forcing}}.
\newblock {Perspectives in Mathematical Logic}. {Springer}, 1998.

\bibitem{Sh:F380}
{Shelah, Saharon}.
\newblock {446 revisited}.

\bibitem{St85}
Jacques Stern.
\newblock Regularity properties of definable sets of reals.
\newblock {\em Annals of Pure and Applied Logic}, 29:289--324, 1985.

\bibitem{Tr77}
John Truss.
\newblock {Sets having calibre $\aleph_1$}.
\newblock In {\em {Logic Colloquium 76}}, volume~87 of {\em Studies in Logic
  and the Foundations of Mathematics}, pages 595--612. North-Holland,
  Amsterdam, 1977.

\end{thebibliography}
\end{document}